# Deadzone-Adapted Disturbance Suppression Control for Global Practical IOS and Zero Asymptotic Gain to Matched Uncertainties

**Iasson Karafyllis[*] and Miroslav Krstic[**]**

[*]Dept. of Mathematics, National Technical University of Athens, Zografou Campus, 15780, Athens, Greece,
email: iasonkar@central.ntua.gr; iasonkaraf@gmail.com

[**]Dept. of Mechanical and Aerospace Eng., University of California, San Diego, La Jolla, CA 92093-0411, U.S.A., email: krstic@ucsd.edu

## Abstract

In this paper we study a special class of systems: time-invariant control systems that satisfy the matching condition for which no bounds for the disturbance and the unknown parameters are known. For this class of systems, we provide a simple, direct, adaptive control scheme that combines three elements: (a) nonlinear damping, (b) single-gain adjustment, and (c) deadzone in the update law. It is the first time that these three tools are combined and the proposed controller is called a Deadzone-Adapted Disturbance Suppression (DADS) Controller. The proposed adaptive control scheme achieves for the first time an attenuation of the plant state to an assignable small level, despite the presence of disturbances and unknown parameters of arbitrary and unknown bounds. Moreover, the DADS Controller prevents gain and state drift regardless of the size of the disturbance and unknown parameter.

## 1. Introduction

Robustness in adaptive control is a major issue that has attracted the attention of many researchers in control theory. Various approaches have been proposed in the literature for time-invariant nonlinear control systems for which no persistence of excitation condition is assumed: nonlinear damping (see [13, 14, 9, 10]), leakage (see [2, 3, 17, 20]), projection methodologies (see [2] and Appendix E in [13]), supervision for direct adaptive schemes (see [1]), dynamic (high) gains or gain adjustment (see [4, 12, 16]) and deadzone in the update law (introduced in the paper [18] and well explained in the book [2]). Every approach has its own advantages and disadvantages and some of the approaches require special assumptions (e.g. knowledge of bounds for the disturbances and/or the unknown parameters).

In the present short note, we focus on a special class of systems: time-invariant systems that satisfy the so-called matching condition (see [13] for the definition of the matching condition) for which no bounds for the disturbance and the unknown parameters are known. For this class of systems, it was recently shown in [10] that there are indirect, adaptive feedback designs which can guarantee KL estimates in the disturbance-free case. Inspired by the idea of dynamic gain used in [4, 19, 12, 16, 15], here we go several steps further. We provide a simple, direct, adaptive control scheme that combines three elements: (a) nonlinear damping (as in [13, 10]), (b) single-gain adjustment (the dynamic feedback has only one state), and (c) deadzone in the update law. It is the first time that these three tools are combined, to the best of the authors knowledge. The proposed adaptive scheme is direct and robustness is not sought by applying advanced identification tools or delays (as in [7, 8, 11]); here no identification is performed (which also explains the simplicity of the proposed scheme). We call the proposed controller a Deadzone-Adapted Disturbance Suppression (DADS) controller.

The advantages of the DADS controller (besides its simplicity) are many. The proposed adaptive control scheme achieves for the first time an attenuation of the plant state to an assignable small level, despite the presence of disturbances and unknown parameters of arbitrary and unknown bounds. Moreover, the DADS controller prevents gain and state drift regardless of the size of the disturbance and unknown parameter. The latter property is a consequence of the Bounded-Input-Bounded-State (BIBS) property and the Input-to-Output Stability (IOS) property that both hold for the closed-loop system. The important (and so rarely achieved) combination of robustness properties that are guaranteed for the closed-loop system in conjunction with the simplicity of the controller justify our focus to a special class of systems (systems with matched uncertainties).

The structure of the paper is as follows. We start with a subsection that provides all the stability notions used in the paper (some of them are modifications of well-known notions presented in [6, 9, 21]). All the main results are stated and discussed in Section 2. Section 3 is devoted to the presentation of certain illustrative examples. Section 4 of the paper contains the proofs of all main results. The concluding remarks of the present work are provided in Section 5.

**Notation and Basic Notions.** Throughout this paper, we adopt the following notation.

* $\mathbb{R}_+ := [0, +\infty)$. For a vector $x \in \mathbb{R}^n$, $|x|$ denotes its Euclidean norm and $x'$ denotes its transpose. We use the notation $x^+$ for the positive part of the real number $x \in \mathbb{R}$, i.e., $x^+ = \max(x, 0)$.

* Let $P \in \mathbb{R}^{n \times n}$ be a symmetric positive definite matrix. By $\lambda_{\min}(P)$ we denote the smallest eigenvalue of $P$.

* Let $D \subseteq \mathbb{R}^n$ be an open set and let $S \subseteq \mathbb{R}^n$ be a set that satisfies $D \subseteq S \subseteq cl(D)$, where $cl(D)$ is the closure of $D$. By $C^0(S; \Omega)$, we denote the class of continuous functions on $S$, which take values in $\Omega \subseteq \mathbb{R}^m$. By $C^k(S; \Omega)$, where $k \geq 1$ is an integer, we denote the class of functions on $S \subseteq \mathbb{R}^n$, which take values in $\Omega \subseteq \mathbb{R}^m$ and have continuous derivatives of order $k$. In other words, the functions of class $C^k(S; \Omega)$ are the functions which have continuous derivatives of order $k$ in $D = int(S)$ that can be continued continuously to all points in $\partial D \cap S$. When $\Omega = \mathbb{R}$ then we write $C^0(S)$ or $C^k(S)$. A function $f \in \bigcap_{k=0}^{\infty} C^k(S; \Omega)$ is called a smooth function.

* By $L^\infty(\mathbb{R}_+)$ we denote the class of essentially bounded, Lebesgue measurable functions $d:\mathbb{R}_+ \to \mathbb{R}^p$. For $d \in L^\infty(\mathbb{R}_+)$ we define $\|d\|_\infty = \sup_{t \geq 0}(|d(t)|)$, where $\sup_{t \geq 0}(|d(t)|)$ is the essential supremum.

* By $K$ we denote the class of increasing continuous functions $a:\mathbb{R}_+ \to \mathbb{R}_+$ with $a(0) = 0$. By $K_\infty$ we denote the class of increasing continuous functions $a:\mathbb{R}_+ \to \mathbb{R}_+$ with $a(0) = 0$ and $\lim_{s \to +\infty}(a(s)) = +\infty$. By $KL$ we denote the set of all continuous functions $\sigma:\mathbb{R}_+ \times \mathbb{R}_+ \to \mathbb{R}_+$ with the properties: (i) for each $t \geq 0$ the mapping $\sigma(\cdot,t)$ is of class $K$; (ii) for each $s \geq 0$, the mapping $\sigma(s,\cdot)$ is non-increasing with $\lim_{t \to +\infty}(\sigma(s,t)) = 0$.

* Let $S \subseteq \mathbb{R}^n$ be a non-empty set with $0 \in S$. We say that a function $V:S \to \mathbb{R}_+$ is positive definite if $V(x) > 0$ for all $x \in S$ with $x \neq 0$ and $V(0) = 0$. We say that a continuous function $V:S \to \mathbb{R}_+$ is radially unbounded if the following property holds: "for every $M > 0$ the set $\{x \in S : V(x) \leq M\}$ is compact". For $V \in C^1(S;\mathbb{R}_+)$ we define $\nabla V(x) = \left(\frac{\partial V}{\partial x_1}(x),\ldots,\frac{\partial V}{\partial x_n}(x)\right)$.

We next recall certain notions of output stability. Let $f:\mathbb{R}^n \to \mathbb{R}^n$ be a locally Lipschitz vector field with $f(0) = 0$ and $h:\mathbb{R}^n \to \mathbb{R}^p$ be a continuous mapping with $h(0) = 0$. Consider the dynamical system

$$\dot{x} = f(x), \, x \in \mathbb{R}^n \tag{0.1}$$

with output

$$Y = h(x) \tag{0.2}$$

We assume that the dynamical system (0.1) is forward complete, i.e., for every $x_0 \in \mathbb{R}^n$ the unique solution $x(t) = \phi(t,x_0)$ of the initial-value problem (0.1) with initial condition $x(0) = x_0$ exists for all $t \geq 0$. We use the notation $Y(t,x_0) = h(\phi(t,x_0))$ for all $t \geq 0$, $x_0 \in \mathbb{R}^n$ and $B_R = \{x \in \mathbb{R}^n : |x| < R\}$ for all $R > 0$. The following properties are standard in the analysis of output stability: see for instance [6, 9, 21]. We say that system (0.1), (0.2) is

i) *Lagrange output stable* if for every $R > 0$ the set $\{|Y(t,x_0)| : x_0 \in B_R, t \geq 0\}$ is bounded.

ii) *Lyapunov output stable* if for every $\varepsilon > 0$ there exists $\delta(\varepsilon) > 0$ such that for all $x_0 \in B_{\delta(\varepsilon)}$, it holds that $|Y(t,x_0)| \leq \varepsilon$ for all $t \geq 0$.

iii) *Globally Asymptotically Output Stable (GAOS)* if system (0.1), (0.2) is Lagrange and Lyapunov output stable and $\lim_{t \to +\infty}(Y(t,x_0)) = 0$ for all $x_0 \in \mathbb{R}^n$.

iv) *Uniformly Globally Asymptotically Output Stable (UGAOS)* if system (0.1), (0.2) is Lagrange and Lyapunov output stable and for every $\varepsilon, R > 0$ there exists $T(\varepsilon,R) > 0$ such that for all $x_0 \in B_R$, it holds that $|Y(t,x_0)| \leq \varepsilon$ for all $t \geq T(\varepsilon,R)$.

It should be noted that (see Theorem 2.2 on page 62 in [6]) UGAOS for system (0.1), (0.2) is equivalent to the existence of a function $\beta \in KL$ such that the following estimate holds for all $x_0 \in \mathbb{R}^n$ and $t \geq 0$:

$$|Y(t,x_0)| \leq \beta(|x_0|,t) \tag{0.3}$$

We say that system (0.1), (0.2) is *practically Uniformly Globally Asymptotically Output Stable (p-UGAOS)* if there exists a function $\beta \in KL$ and a constant $\alpha > 0$ such that the following estimate holds for all $x_0 \in \mathbb{R}^n$ and $t \geq 0$:

$$|Y(t,x_0)| \leq \beta(|x_0|,t) + \alpha \tag{0.4}$$

When $h(x) = x$ then the word "output" in the above properties is omitted (e.g., Lagrange stability, Lyapunov stability, GAS, UGAS, p-UGAS).

Let $f: \mathbb{R}^n \times \mathbb{R}^p \to \mathbb{R}^n$ be a locally Lipschitz vector field with $f(0,0) = 0$. Consider the control system

$$\dot{x} = f(x,d), \ x \in \mathbb{R}^n, \ d \in \mathbb{R}^p \tag{0.5}$$

We assume that system (0.5) is forward complete, i.e., for every $x_0 \in \mathbb{R}^n$ and for every Lebesgue measurable and locally essentially bounded input $d: \mathbb{R}_+ \to \mathbb{R}^p$ the unique solution $x(t) = \phi(t, x_0; d)$ of the initial-value problem (0.5) with initial condition $x(0) = x_0$ corresponding to input $d: \mathbb{R}_+ \to \mathbb{R}^p$ exists for all $t \geq 0$. We use the notation $Y(t, x_0; d) = h(\phi(t, x_0; d))$ for all $t \geq 0$, $x_0 \in \mathbb{R}^n$ and for every Lebesgue measurable and locally essentially bounded input $d: \mathbb{R}_+ \to \mathbb{R}^p$.

We say that system (0.5), (0.2) is *Input-to-Output Stable (IOS)* if there exist functions $\beta \in KL$, $\gamma \in K$ such that the following estimate holds for all $x_0 \in \mathbb{R}^n$, $t \geq 0$ and for every $d \in L^\infty(\mathbb{R}_+)$:

$$|Y(t,x_0;d)| \leq \beta(|x_0|,t) + \gamma(\|d\|_\infty) \tag{0.6}$$

We say that system (0.5), (0.2) is *practically Input-to-Output Stable (p-IOS)* if there exist functions $\beta \in KL$, $\gamma \in K$ and a constant $\alpha > 0$ such that the following estimate holds for all $x_0 \in \mathbb{R}^n$, $t \geq 0$ and for every $d \in L^\infty(\mathbb{R}_+)$:

$$|Y(t,x_0;d)| \leq \beta(|x_0|,t) + \gamma(\|d\|_\infty) + \alpha \tag{0.7}$$

We say that system (0.5), (0.2) satisfies the *practical Output Asymptotic Gain Property (p-OAG)* if there exists a non-decreasing continuous function $\gamma: \mathbb{R}_+ \to \mathbb{R}_+$ with $\gamma(0) = 0$ and a constant $\alpha > 0$ such that the following estimate holds for all $x_0 \in \mathbb{R}^n$, $t \geq 0$ and for every $d \in L^\infty(\mathbb{R}_+)$:

$$\limsup_{t \to +\infty} (|Y(t,x_0;d)|) \leq \gamma(\|d\|_\infty) + \alpha \tag{0.8}$$

When $\gamma \equiv 0$ we say that system (0.5), (0.2) satisfies the *zero practical Output Asymptotic Gain property (zero p-OAG)*.

When $h(x) = x$ then the word "output" in the above properties is either replaced by the word "state" (e.g., ISS, p-ISS) or is omitted (e.g., p-AG, zero p-AG).

## 2. Main Results

In this work we study nonlinear control systems of the form

$$\dot{x} = f(x) + g(x)\left(u + \varphi'(x)\theta + a(x)d\right) \quad (2.1)$$
$$x \in \mathbb{R}^n, u \in \mathbb{R}, d \in \mathbb{R}, \theta \in \mathbb{R}^p$$

where $f, g : \mathbb{R}^n \to \mathbb{R}^n$, $\varphi : \mathbb{R}^n \to \mathbb{R}^p$, $a : \mathbb{R}^n \to \mathbb{R}$ are smooth mappings with $f(0) = 0$, $\varphi(0) = 0$, $x \in \mathbb{R}^n$ is the plant state, $\theta \in \mathbb{R}^p$ is the vector of constant and unknown parameters, $u \in \mathbb{R}$ is the control input and $d \in \mathbb{R}$ is the disturbance. Systems of the form (2.1) are systems that satisfy the so-called matching condition, i.e., the effect of both $\theta \in \mathbb{R}^p$ and $d \in \mathbb{R}$ can be cancelled by the control input $u \in \mathbb{R}$ if they are known. We assume next that $d \in L^\infty(\mathbb{R}_+)$ but we assume no bound for the parameters $\theta \in \mathbb{R}^p$ and the disturbance $d \in \mathbb{R}$.

**Assumption (A):** *There exist smooth mappings* $V, Q : \mathbb{R}^n \to \mathbb{R}_+$, $k : \mathbb{R}^n \to \mathbb{R}$ *with* $k(0) = 0$, $V, Q$ *being positive definite and radially unbounded such that the following inequality holds for all* $x \in \mathbb{R}^n$:

$$\nabla V(x)\left(f(x) + g(x)k(x)\right) \leq -Q(x) \quad (2.2)$$

Assumption 1 guarantees the existence of global feedback stabilizer when $\theta = 0$ and $d = 0$. However, due to the existence of $\theta \in \mathbb{R}^p$ and $d \in \mathbb{R}$ we cannot apply the feedback law $u = k(x)$ without some modification. To this purpose we need an additional assumption.

**Assumption (B):** *There exists a smooth mapping* $\mu : \mathbb{R}^n \to [1, +\infty)$ *such that the following inequality holds for all* $x \in \mathbb{R}^n$:

$$|\varphi(x)|^2 \leq \mu(x)Q(x) \quad (2.3)$$

Assuming that assumptions (A) and (B) hold, we can have the following modifications.

<u>1) Robust Control modification.</u> Assuming that $|\theta| \leq \rho$ where $\rho > 0$ is a constant, the feedback law

$$u = k(x) - c\left(\rho^2 \mu(x) + a^2(x)\right)\nabla V(x)g(x) \quad (2.4)$$

where $c > 1/4$ is a parameter of the controller, guarantees the following inequality for all $x \in \mathbb{R}^n$ and $d \in \mathbb{R}$

$$\nabla V(x)\left(f(x) + g(x)\left(u + \varphi'(x)\theta + a(x)d\right)\right) \leq -\left(1 - \frac{1}{4c}\right)Q(x) + \frac{d^2}{4c} \quad (2.5)$$

Inequality (2.5) is shown by using (2.1), (2.2) and Young inequalities. Inequality (2.5) shows that $V$ is an ISS-Lyapunov function for the closed-loop system (2.1) with (2.4) and the ISS property from the disturbance $d \in L^\infty(\mathbb{R}_+)$ holds. The state is eventually led to a ball around the origin and the radius of the ball is an increasing function of $\frac{1}{4c}\|d\|_\infty^2$.

The disadvantage of the feedback law (2.4) is the fact that it works only when $|\theta| \leq \rho$ and (since $\theta \in \mathbb{R}^p$ is unknown) we cannot guarantee the validity of the inequality $|\theta| \leq \rho$.

2) Adaptive Control modification. In this case, by using the so-called $\sigma$–modification, the dynamic feedback law

$$u = k(x) - \varphi'(x)\hat{\theta} - ca^2(x)\nabla V(x)g(x)$$
$$\frac{d\hat{\theta}}{dt} = \nabla V(x)g(x)\Gamma\varphi(x) - \sigma\hat{\theta} \tag{2.6}$$

where $c, \sigma > 0$ are parameters of the controller and $\Gamma \in \mathbb{R}^{p \times p}$ is a symmetric positive definite matrix guarantees the p-ISS property from the disturbance $d \in L^\infty(\mathbb{R}_+)$ for the closed-loop system (2.1) with (2.6) (with state $(x, \hat{\theta} - \theta) \in \mathbb{R}^n \times \mathbb{R}^p$). Indeed, this is shown by using the ISS-Lyapunov function

$$W(x, \hat{\theta}) = V(x) + \frac{1}{2}(\hat{\theta} - \theta)'\Gamma^{-1}(\hat{\theta} - \theta) \tag{2.7}$$

which satisfies the following inequality for all $x \in \mathbb{R}^n$, $d \in \mathbb{R}$, $\hat{\theta} \in \mathbb{R}^p$ and $\lambda \in (0,1)$

$$\dot{W} \leq -Q(x) - \sigma(1-\lambda)(\hat{\theta} - \theta)'\Gamma^{-1}(\hat{\theta} - \theta) + \frac{\sigma}{4\lambda}\theta\Gamma^{-1}\theta + \frac{d^2}{4c} \tag{2.8}$$

where $\dot{W}$ is the derivative of $W(x, \hat{\theta})$ along the trajectories of the closed-loop system (2.1) with (2.6). Inequality (2.8) shows that the radius of the residual set is an increasing function of $\frac{\sigma}{4\lambda}\theta\Gamma^{-1}\theta + \frac{1}{4c}\|d\|_\infty^2$. This is the disadvantage of the approach: the state is eventually led (even in the disturbance-free case) to a ball around the origin of unknown radius.

In order to present a third option, we need an additional assumption.

**Assumption (C):** *There exist constants $\delta, \eta > 0$ such that the following inequality holds for all $x \in \mathbb{R}^n$ with $|x| \leq \delta$:*

$$Q(x) \geq \eta V(x) \tag{2.9}$$

Assumption (C) allows us to suggest a different option for the robust stabilization of (2.1). Consider the dynamic feedback law:

$$u = k(x) - c(1 + \kappa(\exp(z)))\left(|\varphi(x)|^2 + \lambda^2\mu(x) + a^2(x)\right)\nabla V(x)g(x) \tag{2.10}$$

$$\dot{z} = \Gamma\exp(-z)(V(x) - r)^+ \quad , \quad z \in \mathbb{R} \tag{2.11}$$

where $r, \Gamma, c > 0$, $\lambda \geq 0$ are parameters of the controller (constants) and $\kappa \in K_\infty$ is a smooth function. We call the controller (2.10), (2.11) a Deadzone-Adapted Disturbance Suppression (DADS) controller. The controller (2.10), (2.11) combines the use of deadzone (in (2.11)) and dynamic nonlinear damping (in (2.10)).

The controller (2.10), (2.11) is clearly an extension of the controller (2.4) with $z$ being adapted by means of the update law (2.11) and by treating both $\theta \in \mathbb{R}^p$ and $d \in \mathbb{R}$ as external disturbances: $\theta \in \mathbb{R}^p$ is a constant in time vanishing perturbation while $d \in \mathbb{R}$ is a possibly time-varying, non-vanishing perturbation. The controller (2.10), (2.11) combines both the adaptation idea and the robust control methodology that leads to controller (2.4). However, notice that $\dot{z}$ becomes zero, i.e., the adaptation stops, when the plant state $x$ enters the region defined by $V(x) \leq r$. This is the effect of the deadzone and it is a characteristic that does not appear in the update law of (2.6). The deadzone prevents the state $z$ to grow without bound in the case where the disturbance is present.

The controller (2.10), (2.11) is simple: only one integrator is being used. The dynamic controller gain increases in order to overcome the effect of both $\theta \in \mathbb{R}^p$ and $d \in \mathbb{R}$.

The following theorem clarifies the performance characteristics that the DADS controller (2.10), (2.11) can guarantee for the closed-loop system.

**Theorem 1:** *Suppose that Assumptions (A), (B) and (C) hold. Let $r, \Gamma, c > 0$, $\lambda \geq 0$ be given constants and let $\kappa \in K_\infty$ be a smooth function. Then there exist functions $\beta \in KL$, $\gamma \in K_\infty$ and a non-decreasing function $R: \mathbb{R}_+ \to \mathbb{R}$, all functions independent of $\theta \in \mathbb{R}^p$, such that for every $(x_0, z_0) \in \mathbb{R}^n \times \mathbb{R}$ and for every $d \in L^\infty(\mathbb{R}_+)$ the unique solution of the initial-value problem (2.1), (2.10), (2.11) with initial condition $(x(0), z(0)) = (x_0, z_0)$ is bounded and satisfies the following estimates for all $t \geq 0$:*

$$V(x(t)) \leq \beta(V(x(0)), t) + \gamma\left(\frac{\|d\|_\infty^2 + \left(\left(|\theta| - \lambda\sqrt{c(1+\kappa(\exp(z(0))))}\right)^+\right)^2}{2c(1+\kappa(\exp(z(0))))}\right) \quad (2.12)$$

$$z(0) \leq z(t) \leq R\left(\|d\|_\infty^2 + |\theta|^2 + V(x(0)) + |z(0)|\right) \quad (2.13)$$

$$\limsup_{t \to +\infty}(V(x(t))) \leq r \quad (2.14)$$

**Remarks on Theorem 1:** **(a)** Inequality (2.12) shows that the p-IOS property from the disturbance $d \in L^\infty(\mathbb{R}_+)$ holds for the closed-loop system (2.1), (2.10), (2.11) with output $Y = x$. Indeed, the following estimate is a direct consequence of (2.12):

$$V(x(t)) \leq \beta(V(x(0)), t) + \gamma\left(\frac{\|d\|_\infty^2}{c}\right) + \gamma\left(\frac{|\theta|^2}{c}\right)$$

The above estimate guarantees the p-IOS property from the disturbance $d \in L^\infty(\mathbb{R}_+)$ for the closed-loop system (2.1), (2.10), (2.11) with output $Y = x$.

**(b)** Inequality (2.14) guarantees that the plant state $x$ is eventually led (even in the case where a disturbance is present) to a ball around the origin of _assignable radius_. When the plant state $x$ is considered to be the output of the closed-loop system (2.1), (2.10), (2.11), inequality (2.14) guarantees the zero p-OAG property $\limsup_{t \to +\infty}(|x(t)|) \leq \alpha$ with a constant $\alpha > 0$ independent of the constant parameter $\theta$. Moreover, the state, i.e., both $x(t)$ and $z(t)$, remain bounded for every initial condition and every disturbance $d \in L^\infty(\mathbb{R}_+)$; this is the Bounded-Input-Bounded-State (BIBS) property.

**(c)** The proof of Theorem 1 provides a crude estimate for the function $R$. Let arbitrary $s \geq 0$ be given and define

$$\bar{\eta}(s) := \frac{3}{4} \inf \left\{ \frac{Q(x)}{V(x)} : x \in \mathbb{R}^n, x \neq 0, V(x) \leq 1 + \beta(s,0) + \gamma\left(\frac{s}{2c}\right) \right\}$$

Notice that Assumption (C) guarantees that $\bar{\eta}(s) > 0$. Then we can estimate $R(s)$ by means of the formula

$$R(s) = \ln\left(1 + \kappa^{-1}\left(\left(\frac{s}{4cr\bar{\eta}(s)} - 1\right)^+\right)\right) + \exp(s) + \frac{\Gamma}{\bar{\eta}(s)}\left(\beta(s,0) + \gamma\left(\frac{s}{2c}\right)\right)$$

but the reader is warned that the above formula provides a crude estimate. However, it should be noted that estimate (2.13) in conjunction with estimate (2.12) allows us to obtain a useful qualitative result: that the solution $(x(t), z(t))$ is _uniformly bounded_ for bounded sets of disturbances $d \in L^\infty(\mathbb{R}_+)$ and bounded sets of initial conditions.

**(d)** Taking into account the above remark, we can say that in the disturbance-free case (i.e., when $d \equiv 0$), the closed-loop system (2.1), (2.10), (2.11) with output $Y = x$:
i) is Lagrange stable,
ii) satisfies the p-UGAOS property,
iii) has a globally attracting set, namely the set $\{(x,z) \in \mathbb{R}^n \times \mathbb{R} : V(x) \leq r\}$.

It is interesting to notice that in the disturbance-free case the closed-loop system (2.1), (2.6) with $\sigma = 0$ and output $Y = x$:
i) is Lagrange stable,
ii) is Lyapunov stable,
iii) satisfies the GAOS property.

However, it should also be noted that the closed-loop system (2.1), (2.6) with $\sigma = 0$ does not guarantee the BIBS property when disturbances are present and that is the reason that $\sigma$-modification (i.e., (2.6) with $\sigma > 0$) is being used when disturbances may be present. The closed-loop system (2.1), (2.6) with $\sigma > 0$ satisfies the p-ISS property when disturbances are present and the p-UGAS property in the disturbance-free case with a residual set that depends on $\theta \in \mathbb{R}^p$.

**(e)** For $\lambda > 0$ the DADS controller (2.10), (2.11) has an additional advantage. If there exists $T > 0$ for which $\kappa(\exp(z(T))) \geq \frac{|\theta|^2 - c\lambda^2}{c\lambda^2}$ then estimate (2.12) and the semigroup property imply that the following (non-practical, IOS-like) estimate holds for the closed-loop system (2.1), (2.10), (2.11) with output $Y = x$ and for all $t \geq T$

$$V(x(t)) \leq \beta\left(V(x(T)), t-T\right) + \gamma\left(\frac{\|d\|_\infty^2}{2c}\right)$$

The above estimate shows that in the disturbance-free case, if there exists $T > 0$ for which $\kappa(\exp(z(T))) \geq \frac{|\theta|^2 - c\lambda^2}{c\lambda^2}$ then the plant state is actually led to zero (non-practical regulation). Therefore, in the disturbance-free case we have that: *either the plant state is led asymptotically to zero or the estimate* $\kappa(\exp(z(t))) < \frac{|\theta|^2 - c\lambda^2}{c\lambda^2}$ *holds for all* $t \geq 0$.

**(f)** Theorem 1 is also valid for systems of the form

$$\dot{x} = f(x) + g(x)u + \Phi(x)\theta + A(x)d \quad (2.15)$$
$$x \in \mathbb{R}^n, u \in \mathbb{R}, d \in \mathbb{R}, \theta \in \mathbb{R}^p$$

where $f, g, A: \mathbb{R}^n \to \mathbb{R}^n$, $\Phi: \mathbb{R}^n \to \mathbb{R}^{n \times p}$ are smooth mappings with $f(0) = 0$, $\Phi(0) = 0$ for which there exist smooth mappings $\varphi: \mathbb{R}^n \to \mathbb{R}^p$, $a: \mathbb{R}^n \to \mathbb{R}$ such that $\nabla V(x)\Phi(x) = \nabla V(x)g(x)\varphi'(x)$, $\nabla V(x)A(x) = a(x)\nabla V(x)g(x)$ for all $x \in \mathbb{R}^n$ and for which assumptions (A), (B), (C) hold. The class of systems of the form (2.15) includes systems that do not necessarily satisfy the matching condition.

**(g)** Proposition 4.4 in [5] guarantees that assumptions (A), (B), (C) hold automatically if there exists a smooth mapping $k: \mathbb{R}^n \to \mathbb{R}$ with $k(0) = 0$, such that the feedback law $u = k(x)$ guarantees global asymptotic stability and local exponential stability for system $\dot{x} = f(x) + g(x)k(x)$. In other words, there exist smooth mappings $V, Q: \mathbb{R}^n \to \mathbb{R}_+$, $\mu: \mathbb{R}^n \to [1, +\infty)$ with $V, Q$ being positive definite and radially unbounded for which assumptions (A), (B), (C) hold for system (2.1).

**(h)** Theorem 1 can be stated for less regular mappings $f, g: \mathbb{R}^n \to \mathbb{R}^n$, $\varphi: \mathbb{R}^n \to \mathbb{R}^p$, $V, Q: \mathbb{R}^n \to \mathbb{R}_+$, $k: \mathbb{R}^n \to \mathbb{R}$ and $a: \mathbb{R}^n \to \mathbb{R}$ but here for simplicity reasons we assume that all mappings are smooth.

Assumptions (A), (B), (C) are automatically satisfied for systems of the form

$$\dot{x} = Ax + B\left(u + \varphi'(x)\theta + a(x)d\right) \quad (2.16)$$
$$x \in \mathbb{R}^n, u \in \mathbb{R}, \theta \in \mathbb{R}^p, d \in \mathbb{R}$$

where $A \in \mathbb{R}^{n \times n}$, $B \in \mathbb{R}^n$ is a given stabilizable pair of matrices and $\varphi: \mathbb{R}^n \to \mathbb{R}^p$, $a: \mathbb{R}^n \to \mathbb{R}$ are given smooth mappings with $\varphi(0) = 0$. Using the methodology of the proof of Theorem 1, we obtain the following result.

**Theorem 2:** *Let $A \in \mathbb{R}^{n \times n}$, $B \in \mathbb{R}^n$ be a given stabilizable pair of matrices and let $\varphi: \mathbb{R}^n \to \mathbb{R}^p$, $a: \mathbb{R}^n \to \mathbb{R}$ be given smooth mappings with $\varphi(0) = 0$. Let $k \in \mathbb{R}^n$ be a vector for which $(A - Bk') \in \mathbb{R}^{n \times n}$ is Hurwitz and let $P \in \mathbb{R}^{n \times n}$ be a symmetric positive definite matrix for which $Q = -(A-Bk')'P - P(A-Bk')$ is a positive definite matrix. Let $\eta > 0$ be a constant and let $\mu: \mathbb{R}^n \to [1, +\infty)$ be a smooth function for which the following inequality holds for all $x \in \mathbb{R}^n$:*

$$x'Qx \geq \eta x'Px + \frac{|\varphi(x)|^2}{4\mu(x)} \tag{2.17}$$

*Let $\varepsilon, \Gamma, c > 0$, $\lambda \geq 0$ be given constants and let $\kappa \in K_\infty$ be a smooth function. Define for all $z \in \mathbb{R}$, $a, b \geq 0$ the function:*

$$\chi(z, a, b) := \kappa^{-1}\left( \min\left\{ s \geq \kappa(\exp(z)) : a^2 + \left(\left(b - \lambda\sqrt{c(1+s)}\right)^+\right)^2 \leq 4\eta c \lambda_{\min}(P)\varepsilon^2(1+s) \right\} \right) \tag{2.18}$$

*Consider the feedback law*

$$u = -k'x - 2c\left(1 + \kappa(\exp(z))\right)\left(|\varphi(x)|^2 + \lambda^2 \mu(x) + a^2(x)\right)B'Px$$
$$\dot{z} = \Gamma \exp(-z)\left(x'Px - \lambda_{\min}(P)\varepsilon^2\right)^+ \tag{2.19}$$

*Then for every $(x_0, z_0) \in \mathbb{R}^n \times \mathbb{R}$ and for every $d \in L^\infty(\mathbb{R}_+)$ the unique solution of the initial-value problem (2.16), (2.19) with initial condition $(x(0), z(0)) = (x_0, z_0)$ is bounded and satisfies the following estimates for all $t \geq 0$*

$$x'(t)Px(t) \leq \exp(-\eta t)x'(0)Px(0) + S \tag{2.20}$$

$$\exp(z(0)) \leq \exp(z(t)) \leq \chi\left(z(0), \|d\|_\infty, |\theta|\right) + \frac{\Gamma}{\eta}\left(x'(0)Px(0) + S\right) \tag{2.21}$$

$$\limsup_{t \to +\infty}\left(|x(t)|\right) \leq \varepsilon \tag{2.22}$$

*where* $S = \dfrac{\|d\|_\infty^2 + \left(\left(|\theta| - \lambda\sqrt{c(1 + \kappa(\exp(z(0))))}\right)^+\right)^2}{4\eta c \left(1 + \kappa(\exp(z(0)))\right)}$.

**Remarks on Theorem 2:** (a) Since definition (2.18) implies the following inequality for all $z \in \mathbb{R}$, $a, b \geq 0$

$$\chi(z, a, b) \leq \kappa^{-1}\left( \max\left( \kappa(\exp(z)), \frac{a^2 + \left(\left(b - \lambda\sqrt{c(1+\kappa(\exp(z)))}\right)^+\right)^2}{4\eta c \lambda_{\min}(P)\varepsilon^2} - 1 \right) \right) \tag{2.23}$$

it follows from (2.21) that the following estimate holds for every solution of (2.16), (2.19):

$$\exp(z(t)) \leq \frac{\Gamma}{\eta} x'(0) P x(0) + \Gamma \frac{\|d\|_\infty^2 + \left(\left(|\theta| - \lambda\sqrt{c\left(1 + \kappa(\exp(z(0)))\right)}\right)^+\right)^2}{4\eta^2 c\left(1 + \kappa(\exp(z(0)))\right)}$$

$$+ \kappa^{-1}\left(\max\left(\kappa(\exp(z(0))), \frac{\|d\|_\infty^2 + \left(\left(|\theta| - \lambda\sqrt{c\left(1 + \kappa(\exp(z(0)))\right)}\right)^+\right)^2}{4\eta c \lambda_{\min}(P)\varepsilon^2} - 1\right)\right) \quad (2.24)$$

**(b)** Notice the difference between estimate (2.20) and estimate (2.12): estimate (2.20) guarantees an exponentially vanishing effect of the initial conditions.

**(c)** It should be noticed that systems of the form (2.16) include the chain of integrators

$$\dot{x}_i = x_{i+1} \quad i = 1,\dots,n-1$$
$$\dot{x}_n = u + a(x)d + \sum_{i=1}^{p} \varphi_k(x)\theta_k \quad (2.25)$$
$$x = (x_1,\dots,x_n)' \in \mathbb{R}^n$$

where $\varphi_k : \mathbb{R}^n \to \mathbb{R}$ ($k=1,\dots,p$) and $a : \mathbb{R}^n \to \mathbb{R}$ are smooth functions with $\varphi_k(0) = 0$ ($k=1,\dots,p$).

## 3. Examples

The first example deals with a system of the form (2.25) and presents a comparison between the DADS controller and the adaptive controllers obtained by $\sigma$-modification.

**Example 1:** Consider the planar system

$$\dot{x}_1 = x_2$$
$$\dot{x}_2 = \theta_1 x_1 + \theta_2 x_2 + \theta_3 x_1^2 + u + d \quad (3.1)$$

System (3.1) is a system of the form (2.23) for which we can apply Theorem 2 with

$$n = 2, p = 3, a(x) \equiv 1$$

$$A = \begin{bmatrix} 0 & 1 \\ 0 & 0 \end{bmatrix}, B = \begin{bmatrix} 0 \\ 1 \end{bmatrix}, P = \frac{1}{2}\begin{bmatrix} 5 & 2 \\ 2 & 1 \end{bmatrix}, k = \begin{bmatrix} 5 \\ 4 \end{bmatrix}, Q = 4P$$

$$\varphi(x) = \begin{bmatrix} x_1 & x_2 & x_1^2 \end{bmatrix}', \lambda_{\min}(P) = \frac{3 - 2\sqrt{2}}{2}$$

Moreover, inequality (2.17) holds with $\eta = 3$ and $\mu(x) = 2 + x_1^2$. Using Theorem 2 and the function $\kappa(\rho) = \rho^2$, we know that for every constants $\varepsilon, \Gamma, c > 0$, $\lambda \geq 0$ the dynamic feedback law

$$u = -5x_1 - 4x_2 - c(1+\exp(2z))\left(\left(1+\lambda^2\right)x_1^2 + x_2^2 + x_1^4 + 1 + 2\lambda^2\right)(x_2 + 2x_1)$$

$$\dot{z} = \frac{\Gamma}{2}\exp(-z)\left(x_1^2 + (x_2 + 2x_1)^2 - \left(3 - 2\sqrt{2}\right)\varepsilon^2\right)^+ \tag{3.2}$$

guarantees that for every $(x_0, z_0) \in \mathbb{R}^2 \times \mathbb{R}$, $d \in L^\infty(\mathbb{R}_+)$ the unique solution of the initial-value problem (3.1), (3.2) with initial condition $(x(0), z(0)) = (x_0, z_0)$ satisfies estimate (2.22) as well as the following estimates for $t \geq 0$:

$$x_1^2(t) + (x_2(t) + 2x_1(t))^2 \leq \exp(-3t)\left(x_1^2(0) + (x_2(0) + 2x_1(0))^2\right)$$

$$+ \frac{\|d\|_\infty^2 + \left(\left(|\theta| - \lambda\sqrt{c(1+\exp(2z(0)))}\right)^+\right)^2}{6c(1+\exp(2z(0)))} \tag{3.3}$$

$$\exp(z(0)) \leq \exp(z(t)) \leq \frac{\Gamma}{3}x'(0)Px(0) + \Gamma \frac{\|d\|_\infty^2 + \left(\left(|\theta| - \lambda\sqrt{c(1+\exp(2z(0)))}\right)^+\right)^2}{36c(1+\exp(2z(0)))}$$

$$+ \sqrt{\max\left[\exp(2z(0)), \frac{\|d\|_\infty^2 + \left(\left(|\theta| - \lambda\sqrt{c(1+\exp(2z(0)))}\right)^+\right)^2}{6c\left(3-2\sqrt{2}\right)\varepsilon^2} - 1\right]} \tag{3.4}$$

We next compare the DADS controller (3.2) with the adaptive controller obtained by $\sigma$-modification using the Lyapunov function

$$V(x, \hat{\theta}) = \frac{1}{2}x_1^2 + \frac{1}{2}(x_2 + 2x_1)^2 + \frac{1}{2\Gamma_1}\left(\hat{\theta}_1 - \theta_1\right)^2 + \frac{1}{2\Gamma_2}\left(\hat{\theta}_2 - \theta_2\right)^2 + \frac{1}{2\Gamma_3}\left(\hat{\theta}_3 - \theta_3\right)^2$$

with constants $\Gamma_i > 0$ ($i = 1, 2, 3$) and the nominal feedback law $k(x) = -k'x = -5x_1 - 4x_2$:

$$u = -5x_1 - 4x_2 - c(x_2 + 2x_1) - \hat{\theta}_1 x_1 - \hat{\theta}_2 x_2 - \hat{\theta}_3 x_1^2$$

$$\frac{d\hat{\theta}_1}{dt} = \Gamma_1 x_1 (x_2 + 2x_1) - \sigma\hat{\theta}_1$$

$$\frac{d\hat{\theta}_2}{dt} = \Gamma_2 x_2 (x_2 + 2x_1) - \sigma\hat{\theta}_2 \tag{3.5}$$

$$\frac{d\hat{\theta}_3}{dt} = \Gamma_3 x_1^2 (x_2 + 2x_1) - \sigma\hat{\theta}_3$$

where $c, \sigma > 0$ are constants.

We apply the controllers (3.2) and (3.5) to system (3.1) with $\theta_1 = 10$, $\theta_2 = \theta_3 = 1$. The parameters for the controllers were

$$c = 1, \ \sigma = \Gamma_1 = \Gamma_2 = \Gamma_3 = 1, \ \Gamma = 20, \ \varepsilon = 0.2, \ \lambda = 1$$

The time evolution of the Euclidean norm of the plant state for the solutions of the closed-loop systems (3.1), (3.2) and (3.1), (3.5) is shown in Fig. 1 (disturbance-free case) and in Fig.3 for the disturbance $d(t) = 3\cos(10t)$. It is clear that in both cases the controller (3.2) achieves to bring the plant state very close to 0 while the controller (3.5) fails to achieve this objective.

The time evolution of $\rho(t) = \exp(z(t))$ for the solution of the closed-loop system (3.1), (3.2) is shown in Fig.3. We found that in both cases (the disturbance-free case and the case where $d(t) = 3\cos(10t)$) the evolution is similar: $\rho(t) = \exp(z(t))$ increases and approaches quickly a limit value (3.97 in the disturbance-free case and 3.84 in the case where $d(t) = 3\cos(10t)$). For $t \geq 2$, $\rho(t)$ remains almost constant with value approximately equal to its limit value.

It should also be noticed that Fig. 2 shows that the DADS controller (3.2) achieves almost perfect disturbance rejection (as expected). This is a consequence of the zero p-OAG property from the input $d$ to the plant state $x$ i.e., property (2.22).  ◁

The second example shows that the use of deadzone in (2.11) is absolutely essential for achieving the zero p-OAG property $\limsup_{t \to +\infty}(|x(t)|) \leq \alpha$ from the input $d$ to the plant state $x$ with a constant $\alpha > 0$ independent of the constant parameter $\theta$, i.e., property (2.14).

**Example 2:** As explained above, the mechanism that guarantees the zero p-OAG property $\limsup_{t \to +\infty}(|x(t)|) \leq \alpha$ from the input $d$ to the plant state $x$ with a constant $\alpha > 0$ independent of the constant parameter $\theta$, i.e., property (2.14), in conjunction with the boundedness of the whole state vector $(x, z)$ even in the presence of disturbances, is the combination of the use of deadzone in (2.11) and dynamic nonlinear damping in (2.10). The use of deadzone in (2.11) is absolutely essential. To understand this point, we compare the controller (2.10), (2.11) with a controller for which (2.11) has been replaced by

$$\frac{d}{dt}(\exp(z)) = \Gamma V(x) - \sigma \exp(z) \tag{3.6}$$

where $\sigma > 0$ is a constant. Such a controller may achieve various objectives such as the p-IOS property from the input $d$ to the plant state $x$ and the boundedness of the whole state vector $(x, z)$ even in the presence of disturbances. However, such a controller cannot achieve the zero p-OAG property $\limsup_{t \to +\infty}(|x(t)|) \leq \alpha$ from the input $d$ to the plant state $x$ with a constant $\alpha > 0$ independent of the constant parameter $\theta$, i.e., property (2.14).

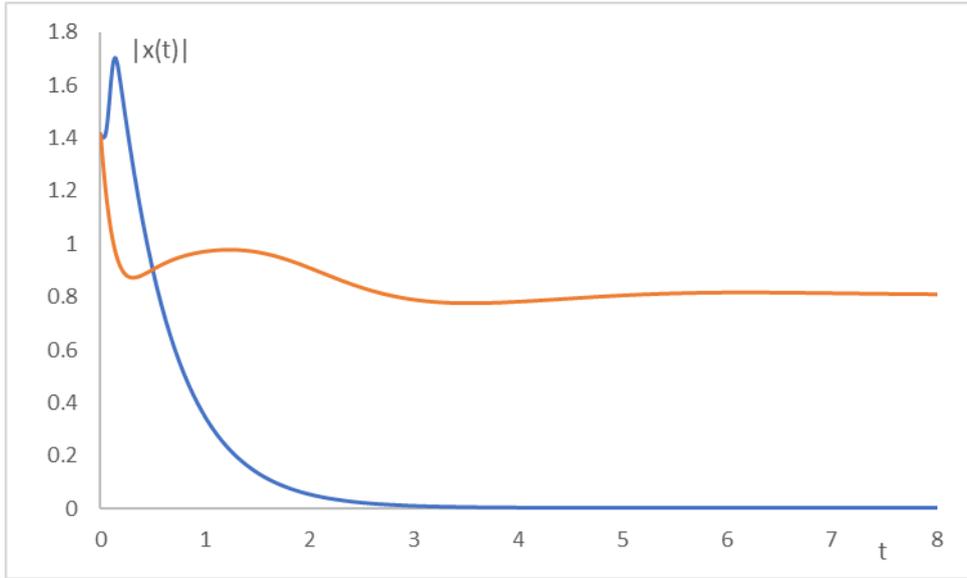

**Fig. 1:** The evolution of the Euclidean norm of the plant state for the solutions of the closed-loop systems (3.1), (3.2) (blue line) and (3.1), (3.5) (red line) with $d \equiv 0$. Initial condition $x_1(0) = -1$, $x_2(0) = 1$, $z(0) = -\ln(10)$, $\hat{\theta}_1(0) = \hat{\theta}_2(0) = \hat{\theta}_3(0) = 0$.

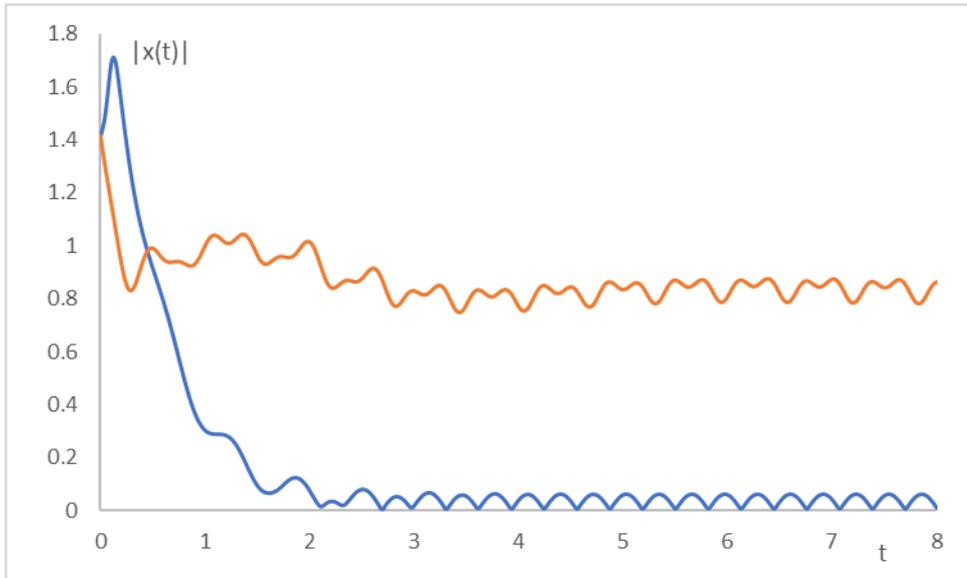

**Fig. 2:** The evolution of the Euclidean norm of the plant state for the solutions of the closed-loop systems (3.1), (3.2) (blue line) and (3.1), (3.5) (red line) with $d(t) = 3\cos(10t)$. Initial condition $x_1(0) = -1$, $x_2(0) = 1$, $z(0) = -\ln(10)$, $\hat{\theta}_1(0) = \hat{\theta}_2(0) = \hat{\theta}_3(0) = 0$.

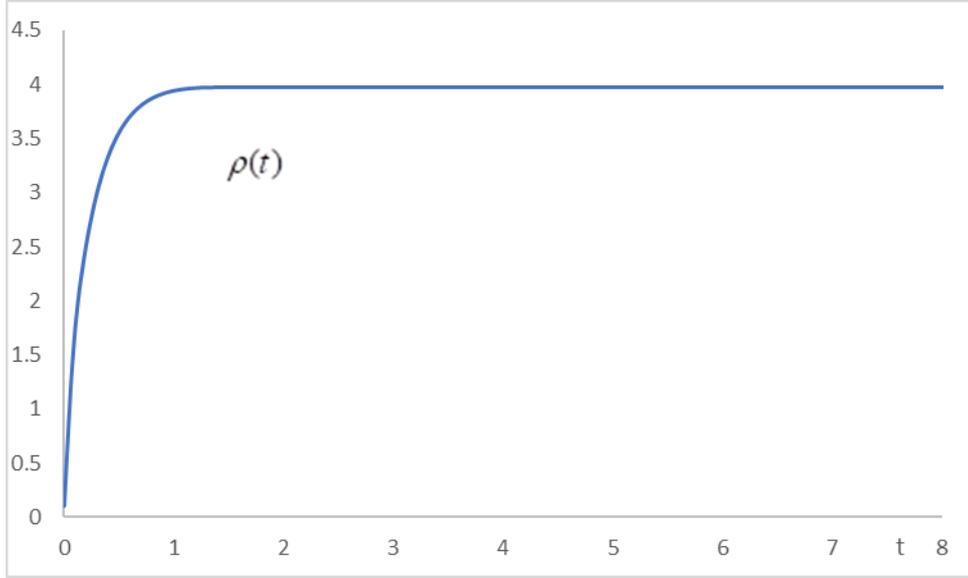

**Fig. 3:** The evolution of $\rho(t) = \exp(z(t))$ for the solution of the closed-loop system (3.1), (3.2) with $d \equiv 0$. Initial condition $x_1(0) = -1$, $x_2(0) = 1$, $z(0) = -\ln(10)$, $\hat{\theta}_1(0) = \hat{\theta}_2(0) = \hat{\theta}_3(0) = 0$.

To see this, let us study the simple scalar example

$$\dot{x} = \theta x + u + d$$
$$x \in \mathbb{R}, u \in \mathbb{R}, \theta \in \mathbb{R}, d \in \mathbb{R} \quad (3.7)$$

A controller based on (2.10), (3.6) and $V(x) = x^2/2$, $\kappa(\rho) = \rho^2$ would give a controller of the form

$$u = -\left(K_1 + K_2\rho^2\right)x - \left(K_3 + K_4\rho^2\right)x^3$$
$$\dot{\rho} = Mx^2 - \sigma\rho \quad (3.8)$$

where $K_1, K_2, K_3, K_4, M, \sigma > 0$ are constants and $\rho = \exp(z)$. The equilibrium points for the closed-loop system (3.7), (3.8) with $d \equiv 0$ are given by the equations

$$\left(K_1 - \theta + K_2\rho^2\right)x + \left(K_3 + K_4\rho^2\right)x^3 = 0$$
$$Mx^2 = \sigma\rho \quad (3.9)$$

It is clear that when $\theta > K_1$ there are non-zero equilibrium points of the form $(x, \rho) = \left(\pm\left(\frac{\sigma\rho}{M}\right)^{1/2}, \rho\right)$ with $\rho > 0$ being the unique positive solution of the equation $\frac{K_4\sigma}{M}\rho^3 + K_2\rho^2 + \frac{K_3\sigma}{M}\rho = \theta - K_1$. Notice that the unique positive solution satisfies the inequality $\rho \geq \min\left(\frac{M(\theta - K_1)}{K_4\sigma + K_2M + K_3\sigma}, \left(\frac{M(\theta - K_1)}{K_4\sigma + K_2M + K_3\sigma}\right)^{1/3}\right)$. Consequently, when $\theta > K_1$ there are solutions of the closed-loop system (3.7), (3.8) with $d \equiv 0$ for which

$$\limsup_{t \to +\infty}(|x(t)|) \geq \left(\frac{\sigma}{M}\right)^{1/2} \min\left(\left(\frac{M(\theta-K_1)}{K_4\sigma+K_2M+K_3\sigma}\right)^{1/2}, \left(\frac{M(\theta-K_1)}{K_4\sigma+K_2M+K_3\sigma}\right)^{1/6}\right)$$

Therefore, it is clear that -no matter how the constants $K_1, K_2, K_3, K_4, M, \sigma > 0$ are selected- property (2.14) cannot hold with a constant $r > 0$ independent of $\theta$.

Therefore, the use of deadzone in (2.11) is absolutely essential for achieving the zero p-OAG property $\limsup_{t \to +\infty}(|x(t)|) \leq \alpha$ from the input $d$ to the plant state $x$ with a constant $\alpha > 0$ independent of the constant parameter $\theta$. ◁

The third example shows that the zero p-OAG property $\limsup_{t \to +\infty}(|x(t)|) \leq \alpha$ from the input $d$ to the plant state $x$ with a constant $\alpha > 0$ independent of the constant parameter $\theta$, i.e., property (2.14), imposes severe structural constraints on the system. In other words, this property cannot be achieved for arbitrary stabilizable systems.

**Example 3:** Consider the planar system
$$\dot{x}_1 = x_2 + d \quad, \quad \dot{x}_2 = u$$
$$x = (x_1, x_2) \in \mathbb{R}^2, u \in \mathbb{R}, d \in \mathbb{R} \quad (3.10)$$

We claim that the zero p-OAG property $\limsup_{t \to +\infty}(|x(t)|) \leq \alpha$ from the input $d$ to the plant state $x$ with a constant $\alpha > 0$ cannot be achieved for system (3.10) (the constant parameter $\theta$ is irrelevant for this example). The proof of this fact is made by contradiction. Suppose that there exists a feedback controller for which the corresponding closed-loop system satisfies the zero p-OAG property $\limsup_{t \to +\infty}(|x(t)|) \leq \alpha$ from the input $d$ to the plant state $x$ with a constant $\alpha > 0$. Consider the solution of (3.10) corresponding to the constant disturbance $d(t) \equiv \alpha + 2$. By assumption, there exists $T > 0$ such that $|x(t)| \leq \alpha + 1$ for $t \geq T$. Using (3.10), we then have for $t \geq T$ a.e.

$$\dot{x}_1(t) = x_2(t) + d(t) \geq \alpha + 2 - |x_2(t)| \geq \alpha + 2 - |x(t)| \geq 1$$

Consequently, we obtain for $t \geq T$ the inequality $\alpha + 1 \geq |x(t)| \geq |x_1(t)| \geq x_1(t) \geq x_1(T) + t - T$, which obviously leads to a contradiction.

It should be noted that the matching condition does not hold for system (3.10). ◁

The final example shows that even if we relax our requirements and try to achieve the p-OAG property $\limsup_{t \to +\infty}(|x(t)|) \leq \alpha(\|d\|_\infty)$ from the input $d$ to the plant state $x$ with a non-decreasing function $\alpha : \mathbb{R}_+ \to (0, +\infty)$ independent of the constant parameter $\theta$, then we still need to impose severe structural constraints on the system.

**Example 4:** Consider the system
$$\dot{x}_1 = -x_1 + d \quad, \quad \dot{x}_2 = \theta x_1 + x_3 \quad, \quad \dot{x}_3 = u$$
$$x = (x_1, x_2, x_3) \in \mathbb{R}^3, u \in \mathbb{R}, \theta \in \mathbb{R}, d \in \mathbb{R} \quad (3.11)$$

We claim that the p-OAG property $\limsup_{t\to+\infty}(|x(t)|) \leq \alpha(\|d\|_\infty)$ from the input $d$ to the plant state $x$ with a non-decreasing function $\alpha : \mathbb{R}_+ \to (0,+\infty)$ independent of the constant parameter $\theta$ cannot be achieved for system (3.11). The proof of this fact is made by contradiction. Suppose that there exists a feedback controller for which the corresponding closed-loop system satisfies the p-OAG property $\limsup_{t\to+\infty}(|x(t)|) \leq \alpha(\|d\|_\infty)$ from the input $d$ to the plant state $x$ with a non-decreasing function $\alpha : \mathbb{R}_+ \to (0,+\infty)$ independent of the constant parameter $\theta$. Consider the solution of (3.11) corresponding to the constant disturbance $d(t) \equiv -1$, constant parameter $\theta = 1 + \alpha(1)$ and initial condition $x_1(0) = 1$. By assumption, there exists $T > 0$ such that $|x(t)| \leq \alpha(1)$ for $t \geq T$. Moreover, $x_1(t) = 1$ for all $t \geq 0$. Using (3.11), we then have for $t \geq T$ a.e.

$$\dot{x}_2(t) = \theta x_1(t) + x_3(t) = \theta + x_3(t)$$
$$\geq \theta - |x_3(t)| \geq \theta - |x(t)| \geq \theta - \alpha(1)$$

Consequently, since $\theta = 1 + \alpha(1)$ we obtain for $t \geq T$ the inequality $\alpha(1) \geq |x(t)| \geq |x_2(t)| \geq x_2(t) \geq x_2(T) + t - T$, which obviously leads to a contradiction.

It should be noted again that the matching condition does not hold for system (3.11). ◁

The last two examples show that the implications of the matching condition are far reaching and should not be neglected by the control theorist or the control practitioner.

## 4. Proofs of Main Results

We next provide the proof of Theorem 1.

**Proof of Theorem 1:** Due to (2.1), (2.10), (2.11) we have for all $(x,z) \in \mathbb{R}^n \times \mathbb{R}$, $d \in \mathbb{R}$:

$$\begin{aligned}\dot{V} &= \nabla V(x)\bigl(f(x) + g(x)k(x)\bigr)\\ &+\nabla V(x)g(x)\varphi'(x)\theta + a(x)\nabla V(x)g(x)d\\ &-c\bigl(1+\kappa(\exp(z))\bigr)\bigl(|\varphi(x)|^2 + \lambda^2\mu(x) + a^2(x)\bigr)\bigl(\nabla V(x)g(x)\bigr)^2\end{aligned} \quad (4.1)$$

where $\dot{V}$ is the derivative of $V(x)$ along the trajectories of the closed-loop system (2.1) with (2.10), (2.11). Using (2.2) and the inequality

$$a(x)\nabla V(x)g(x)d \leq c\bigl(1+\kappa(\exp(z))\bigr)a^2(x)\bigl(\nabla V(x)g(x)\bigr)^2 + \frac{d^2}{4c\bigl(1+\kappa(\exp(z))\bigr)}$$

we get from (4.1) for all $(x,z) \in \mathbb{R}^n \times \mathbb{R}$, $d \in \mathbb{R}$:

$$\dot{V} \le -Q(x) + |\nabla V(x)g(x)||\varphi(x)||\theta| + \frac{d^2}{4c(1+\kappa(\exp(z)))}$$
$$-c(1+\kappa(\exp(z)))\left(|\varphi(x)|^2 + \lambda^2\mu(x)\right)(\nabla V(x)g(x))^2 \quad (4.2)$$

Inequality (4.2) and the fact that $|\theta| \le \left(|\theta| - \lambda\sqrt{c(1+\kappa(\exp(z)))}\right)^+ + \lambda\sqrt{c(1+\kappa(\exp(z)))}$ gives for all $(x,z) \in \mathbb{R}^n \times \mathbb{R}$, $d \in \mathbb{R}$:

$$\dot{V} \le -Q(x) + \left(|\theta| - \lambda\sqrt{c(1+\kappa(\exp(z)))}\right)^+ |\nabla V(x)g(x)||\varphi(x)|$$
$$+ \lambda\sqrt{c(1+\kappa(\exp(z)))} |\nabla V(x)g(x)||\varphi(x)| + \frac{d^2}{4c(1+\kappa(\exp(z)))} \quad (4.3)$$
$$-c(1+\kappa(\exp(z)))\left(|\varphi(x)|^2 + \lambda^2\mu(x)\right)(\nabla V(x)g(x))^2$$

Using the inequalities

$$\lambda\sqrt{c(1+\kappa(\exp(z)))}|\nabla V(x)g(x)||\varphi(x)|$$
$$\le \lambda^2 c(1+\kappa(\exp(z)))\mu(x)(\nabla V(x)g(x))^2 + \frac{1}{4\mu(x)}|\varphi(x)|^2$$

$$\left(|\theta| - \lambda\sqrt{c(1+\kappa(\exp(z)))}\right)^+ |\nabla V(x)g(x)||\varphi(x)|$$
$$\le c(1+\kappa(\exp(z)))|\varphi(x)|^2(\nabla V(x)g(x))^2 + \frac{\left(\left(|\theta| - \lambda\sqrt{c(1+\kappa(\exp(z)))}\right)^+\right)^2}{4c(1+\kappa(\exp(z)))}$$

we get from (4.3) for all $(x,z) \in \mathbb{R}^n \times \mathbb{R}$, $d \in \mathbb{R}$:

$$\dot{V} \le -Q(x) + \frac{1}{4\mu(x)}|\varphi(x)|^2 + \frac{d^2 + \left(\left(|\theta| - \lambda\sqrt{c(1+\kappa(\exp(z)))}\right)^+\right)^2}{4c(1+\kappa(\exp(z)))} \quad (4.4)$$

Inequality (4.4) in conjunction with (2.3) gives for all $(x,z) \in \mathbb{R}^n \times \mathbb{R}$, $d \in \mathbb{R}$:

$$\dot{V} \le -\frac{3}{4}Q(x) + \frac{d^2 + \left(\left(|\theta| - \lambda\sqrt{c(1+\kappa(\exp(z)))}\right)^+\right)^2}{4c(1+\kappa(\exp(z)))} \quad (4.5)$$

Since $V, Q$ are smooth, positive definite and radially unbounded functions, Proposition 2.2 on page 107 in [6] guarantees the existence of a function $\zeta \in K_\infty$ such that the following inequality holds for all $x \in \mathbb{R}^n$:

$$Q(x) \geq \zeta(V(x)) \tag{4.6}$$

Combining (4.5) and (4.6) we get for all $(x,z) \in \mathbb{R}^n \times \mathbb{R}$, $d \in \mathbb{R}$:

$$\dot{V} \leq -\frac{3}{4}\zeta(V(x)) + \frac{d^2 + \left(\left(|\theta| - \lambda\sqrt{c(1+\kappa(\exp(z)))}\right)^+\right)^2}{4c(1+\kappa(\exp(z)))} \tag{4.7}$$

Lemma 2.14 on page 82 in [6] guarantees the existence of $\beta \in KL$ for which the following property holds: for every absolutely continuous function $y:[0,T] \to \mathbb{R}_+$ with $T > 0$ and for every constant $P \geq 0$ for which $\dot{y}(t) \leq -\frac{1}{4}\zeta(y(t))$ holds for almost all $t \in [0,T]$ with $y(t) \geq P$, the following estimate holds for all $t \in [0,T]$:

$$y(t) \leq \max\left(\beta(y(0),t), \beta(P,0)\right) \tag{4.8}$$

Let arbitrary $(x_0, z_0) \in \mathbb{R}^n \times \mathbb{R}$ and $d \in L^\infty(\mathbb{R}_+)$ be given. Consider the unique solution of the initial-value problem (2.1), (2.10), (2.11) with initial condition $(x(0), z(0)) = (x_0, z_0)$. The solution is defined on $[0, t_{max})$, where $t_{max} \in (0, +\infty]$ is the maximal existence time of the solution. Due to (2.11) (which implies that $\dot{z} \geq 0$) it follows that the left inequality (2.13) holds for all $t \in [0, t_{max})$.

The mapping $y(t) = V(x(t))$ is absolutely continuous on every closed interval in $[0, t_{max})$. Due to the left inequality (2.13) and (4.7) we conclude that the differential inequality $\dot{y}(t) \leq -\frac{1}{4}\zeta(y(t))$ holds for almost all $t \in [0, t_{max})$ with $y(t) \geq P = \zeta^{-1}\left(\frac{\|d\|_\infty^2 + \left(\left(|\theta| - \lambda\sqrt{c(1+\kappa(\exp(z(0))))}\right)^+\right)^2}{2c(1+\kappa(\exp(z(0))))}\right)$. It follows from (4.8) that (2.12) holds for all $t \in [0, t_{max})$ with $\gamma(s) := \beta(\zeta^{-1}(s), 0)$ for $s \geq 0$.

Estimate (2.12) implies the estimate

$$V(x(t)) \leq \bar{P} = \beta(V(x(0)), 0) + \gamma\left(\frac{\|d\|_\infty^2 + |\theta|^2}{2c(1+\kappa(\exp(z(0))))}\right) \tag{4.9}$$

for all $t \in [0, t_{max})$. Estimate (4.9) and the fact that $V$ is a radially unbounded function implies that $x(t)$ is bounded on $[0, t_{max})$. Consequently, in order to show that $t_{max} = +\infty$ it suffices to show that $z(t)$ is bounded from above on $[0, t_{max})$. We next show that the right inequality (2.13) holds with

$$R(s) := \ln\left(1 + \max\left(\kappa^{-1}\left(\left(\frac{s}{4cr\bar{\eta}(s)} - 1\right)^+\right), \exp(s)\right)\right) + \frac{\Gamma}{\bar{\eta}(s)}\left(\beta(s,0) + \gamma\left(\frac{s}{2c}\right)\right), \text{ for } s \geq 0 \quad (4.10)$$

where

$$\bar{\eta}(s) := \frac{3}{4}\inf\left\{\frac{Q(x)}{V(x)} : x \in \mathbb{R}^n, x \neq 0, V(x) \leq 1 + \beta(s,0) + \gamma\left(\frac{s}{2c}\right)\right\}, \text{ for } s \geq 0 \quad (4.11)$$

Assumption (C) guarantees that $\bar{\eta}(s) > 0$ for all $s \geq 0$. Notice that definition (4.11) and estimate (4.9) implies that for all $t \in [0, t_{\max})$ it holds that

$$Q(x(t)) \geq \frac{4}{3}\bar{\eta}(s)V(x(t)) \text{ with } s = \|d\|_\infty^2 + |\theta|^2 + V(x(0)) + |z(0)| \quad (4.12)$$

In order to prove the right inequality (2.13) for all $t \in [0, t_{\max})$, we distinguish the following cases.

Case 1: $\exp(z(t)) \leq 1 + \max\left(\kappa^{-1}\left(\left(\frac{s}{4cr\bar{\eta}(s)} - 1\right)^+\right), \exp(s)\right)$ for all $t \in [0, t_{\max})$ with $s = \|d\|_\infty^2 + |\theta|^2 + V(x(0)) + |z(0)|$. Since definition (4.10) implies

$$R(s) \geq \ln\left(1 + \max\left(\kappa^{-1}\left(\left(\frac{s}{4cr\bar{\eta}(s)} - 1\right)^+\right), \exp(s)\right)\right)$$

it follows that the right inequality (2.13) holds for all $t \in [0, t_{\max})$ in this case.

Case 2: There exists $\tau \in [0, t_{\max})$ such that

$$\exp(z(\tau)) > 1 + \max\left(\kappa^{-1}\left(\left(\frac{s}{4cr\bar{\eta}(s)} - 1\right)^+\right), \exp(s)\right).$$

with $s = \|d\|_\infty^2 + |\theta|^2 + V(x(0)) + |z(0)|$. Since

$$\exp(z(0)) < 1 + \max\left(\kappa^{-1}\left(\left(\frac{\|d\|_\infty^2 + |\theta|^2}{4cr\bar{\eta}(s)} - 1\right)^+\right), \exp(z(0))\right)$$

$$\leq 1 + \max\left(\kappa^{-1}\left(\left(\frac{s}{4cr\bar{\eta}(s)} - 1\right)^+\right), \exp(s)\right)$$

it follows that $\tau \neq 0$. Since $z(t)$ is non-decreasing (recall that (2.11) implies that $\dot{z} \geq 0$), there exists $T \in (0, \tau)$ such that $\exp(z(T)) = 1 + \max\left(\kappa^{-1}\left(\left(\frac{s}{4cr\bar{\eta}(s)} - 1\right)^{+}\right), \exp(s)\right)$ and

$$\exp(z(t)) \geq 1 + \max\left(\kappa^{-1}\left(\left(\frac{s}{4cr\bar{\eta}(s)} - 1\right)^{+}\right), \exp(s)\right)$$

for $t \in [T, t_{max})$. Thus, using the fact that $s = \|d\|_{\infty}^2 + |\theta|^2 + V(x(0)) + |z(0)|$, the following inequalities hold for all $t \in [T, t_{max})$

$$\kappa(\exp(z(t))) \geq \left(\frac{s}{4cr\bar{\eta}(s)} - 1\right)^{+} \geq \frac{\|d\|_{\infty}^2 + |\theta|^2}{4cr\bar{\eta}(s)} - 1 \quad (4.13)$$

Consequently, we get from (4.5) and (4.12), (4.13) for $t \in [T, t_{max})$ a.e.:

$$\frac{d}{dt}V(x(t)) \leq -\frac{3}{4}Q(x(t)) + \frac{\|d\|_{\infty}^2 + \left(\left(|\theta| - \lambda\sqrt{c(1+\kappa(\exp(z(t))))}\right)^{+}\right)^2}{4c(1+\kappa(\exp(z(t))))}$$

$$\leq -\frac{3}{4}Q(x(t)) + \frac{\|d\|_{\infty}^2 + |\theta|^2}{4c(1+\kappa(\exp(z(t))))} \leq -\bar{\eta}(s)V(x(t)) + \bar{\eta}(s)r \quad (4.14)$$

Integrating the differential inequality (4.14) we obtain the following estimate for all $t \in [T, t_{max})$:

$$V(x(t)) \leq \exp(-\bar{\eta}(s)(t-T))V(x(T)) + r \quad (4.15)$$

Therefore, we obtain from (2.11) and (4.15) for all $t \in [T, t_{max})$:

$$\frac{d}{dt}(\exp(z(t))) \leq \Gamma \exp(-\bar{\eta}(s)(t-T))V(x(T)) \quad (4.16)$$

Integrating the differential inequality (4.16) we obtain the following estimate for all $t \in [T, t_{max})$:

$$\exp(z(t)) \leq \exp(z(T)) + \frac{\Gamma}{\bar{\eta}}V(x(T)) \quad (4.17)$$

Since $\exp(z(T)) = 1 + \max\left(\kappa^{-1}\left(\left(\frac{s}{4cr\bar{\eta}(s)} - 1\right)^+\right), \exp(s)\right)$, we obtain from (4.17), (4.9) and (4.10) that the right inequality (2.13) holds for all $t \in [T, t_{max})$. Since $z(t)$ is non-decreasing, it follows that the right inequality (2.13) holds for all $t \in [0, t_{max})$.

Consequently, we have proved that (2.13) holds for all $t \in [0, t_{max})$.

By virtue of (2.13), $z(t)$ is bounded on $[0, t_{max})$. Hence, it holds that $t_{max} = +\infty$. Moreover, the solution is bounded for all $t \geq 0$.

Finally, we show inequality (2.14). Since $z(t)$ is non-decreasing and bounded from above the limit $\lim_{t \to +\infty}(\exp(z(t)))$ exists and is finite. Moreover, boundedness of $\dot{x}(t)$ (which follows from (2.1), (2.10), (2.11) and boundedness of the solution) implies that $\frac{d}{dt}(\exp(z(t))) = \Gamma(V(x(t)) - r)^+$ is uniformly continuous. Using Barbălat's Lemma and (2.11) we conclude that $\lim_{t \to +\infty}\left(\frac{d}{dt}(\exp(z(t)))\right) = \lim_{t \to +\infty}\left((V(x(t)) - r)^+\right) = 0$, from which we obtain inequality (2.14). The proof is complete. ◁

We finish this section by providing the proof of Theorem 2.

**Proof of Theorem 2:** Using the functions $V(x) = x'Px$, $k(x) = -k'x$, $Q(x) = x'Qx$, (2.22), (2.17) the fact that $Q = -(A - Bk')'P - P(A - Bk')$ and proceeding exactly as in the proof of Theorem 1, we establish the following inequality for all $(x, z) \in \mathbb{R}^n \times \mathbb{R}$, $d \in \mathbb{R}$:

$$\dot{V} \leq -\eta V(x) + \frac{d^2 + \left(\left(|\theta| - \lambda\sqrt{c(1 + \kappa(\exp(z)))}\right)^+\right)^2}{4c(1 + \kappa(\exp(z)))} \tag{4.18}$$

where $\dot{V}$ is the derivative of $V(x) = x'Px$ along the trajectories of the closed-loop system (2.16) with (2.19).

Let arbitrary $(x_0, z_0) \in \mathbb{R}^n \times \mathbb{R}$ and $d \in L^\infty(\mathbb{R}_+)$ be given. Consider the unique solution of the initial-value problem (2.16), (2.19) with initial condition $(x(0), z(0)) = (x_0, z_0)$. The solution is defined on $[0, t_{max})$, where $t_{max} \in (0, +\infty]$ is the maximal existence time of the solution. Due to (2.19) (which implies that $\dot{z} \geq 0$) it follows that the left inequality (2.21) holds for all $t \in [0, t_{max})$.

The mapping $y(t) = V(x(t)) = x'(t)Px(t)$ is absolutely continuous on every closed interval in $[0, t_{max})$. Due to the left inequality (2.21) and (4.18), we conclude that the differential inequality

$$\dot{y}(t) \leq -\eta\, y(t) + \frac{\|d\|_\infty^2 + |\theta|^2}{4c\left(1+\kappa\left(\exp(z(0))\right)\right)} \tag{4.19}$$

holds for almost all $t \in [0, t_{\max})$. It follows from (4.19) that (2.20) holds for all $t \in [0, t_{\max})$.

Consequently, in order to show that $t_{\max} = +\infty$ it suffices to show that $z(t)$ is bounded from above on $[0, t_{\max})$. Let arbitrary constant $\vartheta > 0$ be given. In order to prove the right inequality (2.21) for all $t \in [0, t_{\max})$, we show first that

$$\exp(z(t)) \leq \vartheta + \chi\left(z(0), \|d\|_\infty, |\theta|\right) + \frac{\Gamma}{\eta} x'(0) P x(0)$$
$$+ \Gamma \frac{\|d\|_\infty^2 + \left(\left(|\theta| - \lambda\sqrt{c\left(1+\kappa\left(\exp(z(0))\right)\right)}\right)^+\right)^2}{4\eta^2 c\left(1+\kappa\left(\exp(z(0))\right)\right)}, \quad \text{for all } t \in [0, t_{\max}) \tag{4.20}$$

We distinguish the following cases.

Case 1: $\exp(z(t)) \leq \vartheta + \chi\left(z(0), \|d\|_\infty, |\theta|\right)$ for all $t \in [0, t_{\max})$. It follows that inequality (4.20) holds in this case.

Case 2: There exists $\tau \in [0, t_{\max})$ such that

$$\exp(z(\tau)) > \vartheta + \chi\left(z(0), \|d\|_\infty, |\theta|\right)$$

Since definition (2.18) implies that $\chi\left(z(0), \|d\|_\infty, |\theta|\right) \geq \exp(z(0))$, it follows that $\tau \neq 0$ and that $\exp(z(0)) < \vartheta + \chi\left(z(0), \|d\|_\infty, |\theta|\right)$. Since $z(t)$ is non-decreasing (recall that (2.19) implies that $\dot{z} \geq 0$), there exists $T \in (0, \tau)$ such that $\exp(z(T)) = \vartheta + \chi\left(z(0), \|d\|_\infty, |\theta|\right)$ and

$$\exp(z(t)) \geq \vartheta + \chi\left(z(0), \|d\|_\infty, |\theta|\right) \tag{4.21}$$

for $t \in [T, t_{\max})$. Definition (2.18) implies that

$$\kappa\left(\chi\left(z(0), \|d\|_\infty, |\theta|\right)\right) \geq \frac{\|d\|_\infty^2 + \left(\left(|\theta| - \lambda\sqrt{c\left(1+\kappa\left(\chi\left(z(0), \|d\|_\infty, |\theta|\right)\right)\right)}\right)^+\right)^2}{4\eta c \lambda_{\min}(P)\varepsilon^2} - 1 \tag{4.22}$$

Combining (4.21) and (4.22) we get for $t \in [T, t_{\max})$

$$1+\kappa(\exp(z(t))) \geq \frac{\|d\|_\infty^2 + \left(\left(|\theta|-\lambda\sqrt{c\left(1+\kappa\left(\chi\left(z(0),\|d\|_\infty,|\theta|\right)\right)\right)}\right)^+\right)^2}{4\eta c\lambda_{\min}(P)\varepsilon^2} \qquad (4.23)$$

$$\geq \frac{\|d\|_\infty^2 + \left(\left(|\theta|-\lambda\sqrt{c\left(1+\kappa(\exp(z(t)))\right)}\right)^+\right)^2}{4\eta c\lambda_{\min}(P)\varepsilon^2}$$

Consequently, we get from (4.18) and (4.23) for $t \in [T, t_{\max})$ a.e.:

$$\frac{d}{dt}V(x(t)) \leq -\eta V(x(t)) + \frac{\|d\|_\infty^2 + \left(\left(|\theta|-\lambda\sqrt{c\left(1+\kappa(\exp(z(t)))\right)}\right)^+\right)^2}{4c\left(1+\kappa(\exp(z(t)))\right)} \qquad (4.24)$$

$$\leq -\eta V(x(t)) + \lambda_{\min}(P)\varepsilon^2 \eta$$

Integrating the differential inequality (4.24) we obtain the following estimate for all $t \in [T, t_{\max})$:

$$V(x(t)) \leq \exp(-\eta(t-T))V(x(T)) + \lambda_{\min}(P)\varepsilon^2 \qquad (4.25)$$

Therefore, we obtain for all $t \in [T, t_{\max})$ from (2.19), (4.25) and the fact that $V(x(t)) = x'(t)Px(t)$:

$$\frac{d}{dt}(\exp(z(t))) \leq \Gamma \exp(-\eta(t-T))V(x(T)) \qquad (4.26)$$

Integrating the differential inequality (4.26) we obtain the following estimate for all $t \in [T, t_{\max})$:

$$\exp(z(t)) \leq \exp(z(T)) + \frac{\Gamma}{\eta}V(x(T)) \qquad (4.27)$$

Since $\exp(z(T)) = \vartheta + \chi(z(0), \|d\|_\infty, |\theta|)$, we obtain from (4.27), (2.20) and the fact that $V(x(t)) = x'(t)Px(t)$ that inequality (4.24) holds for all $t \in [T, t_{\max})$. Since $z(t)$ is non-decreasing, it follows that the right inequality (4.20) holds for all $t \in [0, t_{\max})$.

Consequently, we have proved that (4.20) holds for all $t \in [0, t_{\max})$ and for arbitrary $\vartheta > 0$. Therefore, (2.21) holds for all $t \in [0, t_{\max})$. Hence, by virtue of (2.21), $z(t)$ is bounded from above on $[0, t_{\max})$. Hence, it holds that $t_{\max} = +\infty$. Moreover, the solution is bounded for all $t \geq 0$.

Finally, we show inequality (2.22). Since $z(t)$ is non-decreasing and bounded from above the limit $\lim_{t \to +\infty}(\exp(z(t)))$ exists and is finite. Moreover, boundedness of $\dot{x}(t)$ (which follows from (2.16), (2.19), and boundedness of the solution) implies that $\frac{d}{dt}(\exp(z(t))) = \Gamma\left(x'(t)Px(t) - \lambda_{\min}(P)\varepsilon^2\right)^+$ is uniformly continuous. Using Barbălat's Lemma and (2.19) we conclude that $\lim_{t\to+\infty}\left(\frac{d}{dt}(\exp(z(t)))\right) = \lim_{t\to+\infty}\left((x'(t)Px(t) - \lambda_{\min}(P)\varepsilon^2)^+\right) = 0$, from which inequality (2.22) follows. The proof is complete. ◁

## 5. Concluding Remarks

As mentioned in the Introduction, the important (and so rarely achieved) combination of robustness properties that are guaranteed by using the DADS controller for the closed-loop system in conjunction with its simplicity justify our focus to a special class of systems (systems with matched uncertainties). However, it is clear that a detailed study is needed for the characterization of the class of nonlinear systems for which these robustness properties can be guaranteed by adaptive control schemes. Moreover, it will be useful to construct more academic examples for which there is no adaptive control scheme that guarantees for the corresponding closed-loop system: (a) the p-IOS property from the disturbance to the plant state, (b) the BIBS property, and (c) the zero p-OAG property $\limsup_{t\to+\infty}(|x(t)|) \leq \alpha$ with a constant $\alpha > 0$ independent of the constant parameters.